\title{Every Minor-Closed Property of Sparse Graphs is Testable}
\author{Itai Benjamini \thanks{The Weizmann Institute and Microsoft Research}
\and Oded Schramm \thanks{Microsoft Research}
\and Asaf Shapira \thanks{Microsoft Research}}

\date{} 
\documentclass[11pt]{article}
\usepackage{amsmath}
\usepackage{amssymb}
\usepackage{amsthm}
\usepackage{amsfonts}
\usepackage{graphicx}

\newif\iffigures\figurestrue
\figuresfalse

\newif\ifhyper\IfFileExists{hyperref.sty}{\hypertrue}{\hyperfalse}

\ifhyper\usepackage{hyperref}
\def\hitem#1#2{\item[\hypertarget{#1}{#2}]\expandafter\gdef\csname LBL#1ITM\endcsname{#2}}
\def\iref#1{\hyperlink{#1}{\csname LBL#1ITM\endcsname}}
\else
\def\hitem#1#2{\item[{#2}]\expandafter\gdef\csname LBL#1ITM\endcsname{#2}}
\def\iref#1{{\csname LBL#1ITM\endcsname}}
\fi

\newif\ifdraft
\drafttrue

\long\def\comment#1{}
\long\def\old#1{}

\numberwithin{equation}{section}
\numberwithin{figure}{section}

\newtheorem{theorem}{Theorem}
\numberwithin{theorem}{section}
\newtheorem{corollary}[theorem]{Corollary}
\newtheorem{lemma}[theorem]{Lemma}
\newtheorem{proposition}[theorem]{Proposition}
\newtheorem{conjecture}[theorem]{Conjecture}
\newtheorem{problem}[theorem]{Problem}

\newtheorem*{defi}{Definition} 

\let\qqed=\qed
\def\QED{\qqed\medskip}
\let\qed=\QED

\newcommand{\N}{\mathbb{N}}

\def\Var{\mathrm{Var}}
\def \eps {\epsilon}
\def \P {{\bf P}}
\def\md{\mid}
\def\Bb#1#2{{\def\md{\bigm| }#1\bigl[#2\bigr]}}
\def\BB#1#2{{\def\md{\Bigm| }#1\Bigl[#2\Bigr]}}
\def\Bs#1#2{{\def\md{\mid}#1[#2]}}
\def\Pb{\Bb\P}
\def\Eb{\Bb\E}

\def\EB{\BB\E}
\def\Ps{\Bs\P}

\def \p {{\partial}}
\def \E {{\bf E}}

\def \proof {{ \medbreak \noindent {\bf Proof.} }}
\def\proofof#1{{ \medbreak \noindent {\bf Proof of #1.} }}

\def\dd{\rho}

\def\noopsort#1{}

\begin{document}
\maketitle

\begin{abstract}
Suppose $G$ is a graph with degrees bounded by $d$, and one needs to
remove more than $\epsilon n$ of its edges in order to make it planar. We show
that in this case the statistics of local neighborhoods around
vertices of $G$ is far from the statistics of local neighborhoods
around vertices of any planar graph $G'$ with the same degree bound.
In fact, a similar result
is proved for any minor-closed property of bounded degree graphs.

As an immediate corollary of the above result we infer that many
well studied graph properties, like being planar, outer-planar,
series-parallel, bounded genus, bounded tree-width and several
others, are testable with a constant number of queries,
where the constant may depend on $\epsilon$ and $d$, but not on 
the graph size. None of
these properties was previously known to be testable even with
$o(n)$ queries.
\end{abstract}


\section{Introduction}\label{intro}

Suppose we are given an $n$-vertex graph of bounded degree and are
asked to decide if it is planar. This problem is well known to be
solvable in time $\Theta(n)$ \cite{HT}. But suppose we are only
asked to distinguish with probability $2/3$ between the case that the
input is planar from the case that an $\epsilon$-fraction of its
edges should be removed in order to make it planar. Can we design a
faster algorithm for this relaxed version of the problem, with
running time $o(n)$? A special case of our main result is that in
this case we can even design an algorithm whose running time is a
{\em constant} that depends only on $\epsilon$ (and the bound on the
degrees) and is independent of the size of the input.

\subsection{Background on property testing}

Before stating our main result, which is a significant
generalization of the above mentioned result on planarity testing,
let us briefly introduce the basic notions in the field of property
testing. The meta problem is the following: given a combinatorial
structure $S$, distinguish if $S$ satisfies some property ${\cal P}$
or if $S$ is $\epsilon$-far from satisfying ${\cal P}$, where $S$ is
said to be $\epsilon$-far from satisfying ${\cal P}$ if an
$\epsilon$-fraction of its representation should be modified in
order to make $S$ satisfy ${\cal P}$. The main goal is to design
randomized algorithms, which look at a very small portion of the
input, and using this information distinguish with high probability
between the above two cases. Such algorithms are called {\em
property testers} or simply {\em testers} for the property ${\cal
P}$. Preferably, a tester should look at a portion of the input
whose size is a function of $\epsilon$ only. Blum, Luby and
Rubinfeld \cite{BLR} were the first to formulate a question of this
type, and the general notion of property testing was first
formulated by Rubinfeld and Sudan \cite{RS}.

The main focus of this paper is testing properties of graphs in the
{\em bounded degree} model, which was introduced by Goldreich and
Ron \cite{GR}. In this model, we fix a degree bound $d$ and represent
graphs using adjacency lists. More precisely, we assume that a graph
$G$ is represented as a function $f_G:[n]\times[d]\mapsto
[n]\cup\{*\}$, where given a vertex $v \in V(G)$ and $1 \leq
i\leq d$ the function $f(v,i)$ returns the $i^{th}$ neighbor of $v$,
in case $v$ has at least $i$ vertices. If $v$ has less than $i$
vertices then $f(v,i)=*$. A graph of bounded degree $d$ is said to
be $\epsilon$-far from satisfying ${\cal P}$ if one needs to execute
at least $\eps d n$ edge operations of deleting or adding an edge
to $G$ in order to turn
it into a graph satisfying ${\cal P}$.
(Since $d$ is treated as a constant in this paper,
the $\epsilon d n$ above is just proportional to $\epsilon n$.)

A testing algorithm (or tester) for graph property\footnote{As
usual, a graph property is simply a family of graphs closed under
graph isomorphism.} ${\cal P}$ and accuracy $\eps$ is a (possibly
randomized) algorithm that distinguishes with probability at least
$2/3$ between graphs satisfying ${\cal P}$ from graph that are
$\epsilon$-far from satisfying it. More precisely, if the input
graph satisfies ${\cal P}$ the algorithm accepts it with probability
at least $2/3$, where the probability is taken over the coin tosses
of the tester. Similarly, if it is $\epsilon$-far from satisfying
${\cal P}$, the algorithm should reject it with probability at least
$2/3$.
The tester is given $n$ as input\footnote{%
Our tester needs $n$ as the input only to be
able to pick a vertex at random from the tested graph.}
and is provided with
access to the function $f_G$ as a black box.
We define the query complexity $q_{\mathcal T}(n)$ of the
tester $\mathcal T$ as the maximal number of $f_G$-calls
that the tester can execute on any graph $G$ with $n$ vertices.
Finally, we define the following notion of
efficient testing:

\begin{defi}[Testable]
A graph property ${\cal P}$ is {\em testable} if
for every $\eps>0$ there is some constant $c(\eps)<\infty$
and a tester $\mathcal T$
that tests $\mathcal P$ with accuracy $\eps$
and satisfies $\sup_n q_{\mathcal T}(n)<c(\eps)$ for every $n$.
In other words, if for any $\epsilon > 0$, there
is a {\em constant} time randomized algorithm that can distinguish
with high probability between graphs satisfying ${\cal P}$ from
those that are $\epsilon$-far from satisfying it.
\end{defi}

\subsection{Minor-closed properties and the main result}

Our main result deals with the testing of minor closed graph
properties. Let us briefly introduce the basic notions in this area,
which is too rich to thoroughly survey here. For more details, see
Chapter 12 of \cite{Dies} and the recent article of Lov\'asz
\cite{L} on the subject. A graph $H$ is said to a {\em minor} of a
graph $G$, if $H$ can be obtained from $G$ using a sequence of
vertex removals, edge removals and edge
contractions\footnote{Contracting an edge connecting vertices $u,v$
is the result of replacing $u$ and $v$ by a new vertex $w$, and
connecting $w$ to all the vertices that were connected to either $u$
or $v$. Loops and multiple edges are removed.}. If $H$ is not a
minor of $G$ then $G$ is said to be {\em $H$-minor free}. The topic
of graph minors is among the most (or perhaps the single most)
studied concepts in graph theory. Our main result in this paper is
accordingly the following:

\begin{theorem}[Main Result]\label{t1} For
every (finite) graph $H$, the property of being $H$-minor free is testable.
More generally, every minor-closed graph property is testable.
\end{theorem}

Recall that a graph property $\mathcal P$ is {\em minor-closed}
if every minor of a graph in $\mathcal P$ is also in $\mathcal P$, or equivalently
if ${\cal P}$ is closed under removal of edges, removal of vertices and contraction
of vertices.

Perhaps the most well-known result in the area of graph minors is
the Kuratowski-Wagner Theorem \cite{Kur,W}, which states that a
graph is planar if and only if it is $K_5$-minor free and
$K_{3,3}$-minor free. This fundamental result raised the natural
question if a similar characterization, using a finite family of
forbidden minors, also holds for embedding graphs in other fixed
surfaces? Observe that a graph remains planar if we remove one of
its edges or vertices and if we contract one of its edge. In fact,
this closure property under these three basic operations holds for
the property of a graph being embeddable in any specific surface.
Thus the property of being embeddable in a specific surface is
minor-closed. In one of the deepest results in graph theory,
Robertson and Seymour proved the so called Graph-Minor Theorem
\cite{RobSey2}, which states that for every minor-closed graph
property ${\cal P}$, there is a finite family of graphs ${\cal
H}_{\cal P}$ such that a graph satisfies ${\cal P}$ if and only if
it is $H$-minor free for all $H \in {\cal H}_{\cal P}$. Note that
this in particular answers the above mentioned problem regarding the
characterization of graphs embeddable in a fixed surface using a
finite number of forbidden minors.

We note that besides the above (mainly) graph theoretic motivation,
the graph minors and minor closed properties have also received a
considerable amount of attention by the Computer-Science community,
because many natural graph properties (e.g., planarity)
are described by excluded minors, and graph problems that
are NP-hard in general can often be solved on excluded minor
families in polynomial time, or at least be approximated
better than on arbitrary graphs.
See \cite{DHK} and its references.

Let us mention some well studied minor-closed graph properties. Of
course, the most well known such property is {\bf Planarity}. A
well-studied variant of planarity is {\bf Outer-planarity}, the
property of being embeddable in the plane in such a way that all
vertices lie on the outer face. Another generalization of planarity
is being embeddable in a surface of genus
at most $k$. Graphs satisfying this property are said to
have {\bf genus $k$}.
The {\bf Tree-width} of a graph is one of the most important
invariants of graphs. This notion, which measures how close a graph
is to being a tree, was introduced by Robertson and Seymour as
part of their proof of the Graph-Minor Theorem; see \cite{Bod}. It
also has numerous applications in the area of fixed-parameter
algorithms; see \cite{DF} for more details. As it turns out, having
bounded tree-width is also a minor-closed graph property. Another
well-known minor-closed property is being {\bf Series-parallel}.
Series parallel graphs are graphs that can be obtained from a single
edge by sequence of parallel extensions (adding an edge parallel to
an edge that already exists) and series extensions (subdividing an
edge by a new node). We conclude with the property of being {\bf
Knotlessly-embeddable}, which is the property of being embeddable in
$\mathbb{R}^3$ in a way that no two cycles are linked and no cycle
is knotted.

By Theorem \ref{t1} we infer that the above mentioned properties are
all testable with a constant number of queries. It is interesting to
note that prior to this work none of the above properties was even
known to be testable with $q_{\mathcal T}(n) = o(n)$
queries\footnote{The conference version of \cite{GR} announced a
tester for planarity, which eventually turned out to be wrong.}.

One important aspect of Theorem \ref{t1}, which we have neglected to
address thus far, is the actual dependence of the query-complexity
of the testers on $\epsilon$. As it does not seem like one can
achieve a query complexity that is sub-exponential in
$1/\epsilon$ using our approach, we
opted in several places not to give explicit bounds. We further
discuss this issue in Section \ref{open}, where we show how can one
derive explicit (but rather large) upper bounds on the query
complexity of the testers as a function of $\epsilon$.

\subsection{Hyper-finiteness and testing monotone hyper-finite properties}

The following notion of hyper-finiteness was defined by
Elek \cite{Elek} (though it is implicit in~\cite{LT1}, for example).

\begin{defi}[Hyper-Finite]
A graph $G=(V,E)$
is $(\delta,k)$-hyper-finite if one can remove
$\delta |V|$ edges from $G$ and obtain a graph with connected
components of size at most $k$.
A collection of graphs $\mathcal G$ is hyper-finite if
for every $\delta>0$ there is some finite $k=k(\delta)$
such that every graph in $\mathcal G$ is $(\delta,k)$-hyper-finite.
\end{defi}

Our main theorem, Theorem~\ref{t1}, will be a corollary of the
following more general result.

\begin{theorem}\label{t.m}
Every monotone hyper-finite graph property is testable.
\end{theorem}

Recall that a graph property $\mathcal P$ is {\em monotone}
if every subgraph of a graph in $\mathcal P$ is also in $\mathcal P$.

For example, consider the property that the number of
vertices at distance $r$ around every vertex is at most
some fixed function $g(r)$ satisfying $g(r)= 2^{o(r)}$ as $r\to\infty$.
Theorem~\ref{t.m} can be used to show that this property is testable
(but this does not seem to follow directly from Theorem~\ref{t1}).
Elek \cite{Elek2} studied property testing within this class of graphs.

\subsection{Comparison to previous results}\label{compare}

Our main result in this paper deals with testing properties of
bounded degree graphs. Let us briefly mention some results on
testing properties of dense graphs, a model that was introduced by
Goldreich, Goldwasser and Ron \cite{GGR}. In this model, a graph $G$
is said to be $\epsilon$-far from satisfying a property ${\cal P}$,
if one needs to add/delete at least $\epsilon n^2$ edges to/from $G$
in order to turn it into a graph satisfying ${\cal P}$.
 The tester can ask an oracle whether a pair of vertices, say
$i$ and $j$, are adjacent in the input graph $G$. It was shown in
\cite{GGR} that every partition problem like $k$-colorability and
Max-Cut is testable in this model. Alon and Shapira \cite{ASHer}
have shown that every hereditary graph property\footnote{A graph
property is hereditary if it is closed under removal of vertices.
Therefore, any minor closed property is hereditary} is testable in
dense graphs. This also gave an (essential) characterization of the
graph properties that are testable with one-sided error\footnote{A
tester has one-sided error if it always accepts graphs satisfying
the property.}. A characterization of the properties that are
testable in dense graphs was obtained by Alon et al.\ \cite{AFNS}.
Note that in this model (as its name suggests) we implicitly assume
that the input graph is dense, because the definition of
$\epsilon$-far is relative to $n^2$. Therefore, some properties are
trivially testable in this model. In particular, minor-closed
properties are trivially testable in this model with $O(1/\epsilon)$
queries and even with one-sided error. This follows from the result
of Kostochka and independently Thomason \cite{K1,K2,T1,T2} that
every finite graph with average degree at least $\Omega(r\sqrt{\log
r})$ contains every graph on $r$ vertices as a minor. Therefore,
every large enough finite graph with $\Omega(n^2)$ edges does not
satisfy a minor-closed property.

As the above mentioned results indicate, testing properties in dense
graphs is relatively well understood. In sharp contrast, there were
no general results on testing properties in the bounded-degree
model. For example, while every hereditary property is testable in
dense graphs, in bounded degree graphs some properties are testable
(e.g.\ being triangle-free \cite{GR}), some require
$\tilde{\Theta}(\sqrt{n})$ queries (e.g.\ being bipartite
\cite{GR,GR99}) and some require $\Theta(n)$ queries (e.g.\
3-colorability \cite{BOT}). Besides the above mentioned results, it
was also shown in \cite{GR} that $k$-connectivity is testable in
bounded degree graphs, and Czumaj and Sohler \cite{CzSo} have
recently shown that a relaxed version of expansion is testable with
$\tilde{\Theta}(\sqrt{n})$ queries. Finally, Czumaj, Sohler and
Shapira \cite{CSS} have recently shown that every hereditary
property is testable if the input graph is guaranteed to be non-expanding. Some of
the arguments in the present paper are motivated by some of the
ideas from \cite{CSS}.

One reason for the fact that we understand testing of dense graphs
more than we understand testing of bounded-degree graphs is that
there are structural results ``describing'' dense graphs, primarily
Szemer\'edi's regularity lemma \cite{Sz}, while there are no similar
results for arbitrary sparse graphs. As is evident from the above
discussion, our main result here is the first to show that a general
(and natural) family of properties are all testable in
bounded-degree graphs.

For more details, see the surveys \cite{ASSurvey,CzSo2,F,Rub,Ron}
and Section \ref{open} where we discuss another model of testing
graph properties.

\subsection{Techniques and overview of the paper}

In the next section, we introduce a metric condition for testability.
Basically, we define a sequence of pseudometrics $\dd_r$ indexed
by an integer $r\in\N_+$, where $\dd_r(G,G')$ measures
the difference between $G$ and $G'$ in the frequency of
isomorphism types of $r$-neighborhoods of vertices.
We show that if for some fixed $r>0$ the $\dd_r$ distance
between two graph families is positive, then these graph
families can be distinguished by a tester.

Next, we state a theorem (Theorem~\ref{t.gap}) that roughly says that
$(\eps,k)$-hyper-finite graphs are far away in some
pseudometric $\dd_R$ from graphs that are
not $(\eps',k)$-hyper-finite, where $\eps'$ depends
on $\eps$ and tends to zero as $\eps\searrow 0$.
This result can be deduced from a
recent result of Schramm \cite{Sc}, concerning
properties of convergent sequences of bounded degree graphs.
However, we provide an alternative self-contained proof
in Section \ref{SecGap}. In contrast with~\cite{Sc},
our proof is finitary and gives an explicit upper bound
on $R$, while the proof in \cite{Sc} did not supply such a bound. We note
that convergent sequences of graphs have been previously used in the
study of property testing of {\em dense} graphs \cite{BCLSSV06}.

How does a tester for a monotone hyperfinite graph property $\mathcal P$
work? A first guess might be the following: as the property is monotone,
we may expect to sample enough vertices such that with high probability the neighborhood
of at least one of them will not satisfy the property, thus establishing that the graph itself
does not satisfy the property. This works for properties like triangle-freeness, but
it is not difficult to see, however (see item 2, in the concluding remarks)
that this approach of aiming for a one-sided error test is bound to fail in our case.
The reason is that a graph can be far from being (say) planar, yet locally be planar.
But still, the local neighborhoods of vertices do tell us something about its global properties.
For example, in \cite{GR}, it was shown that the local density around vertices can be used to test the property
of being cycle-free. By the above discussion we must look for another property that the graph
must satisfy if it satisfies a minor closed property. As it turns out, the right property to look
for is hyper-finiteness.

We now proceed with an informal and not very precise description
of how and why our tester works. Given an input graph $G$, the tester first tests for hyper-finiteness.
The test for hyper-finiteness proceeds by picking a bounded
number of vertices at random and exploring the $R$-neighborhood of
each one (where the number of vertices chosen and $R$ depend on $\eps$).
It then checks if the observed frequency of the isomorphism types
of these neighborhoods looks approximately like some $(\eta,k)$-hyper-finite
graph, with appropriate values for $\eta$ and $k$.
We do not have an explicit characterization of the approximate frequencies
occurring in appropriately hyper-finite graphs, but since only an approximation
is needed (with known accuracy), there is a finite table listing the occurring
frequencies, which can be part of the algorithm. The
point here is that the table does not depend on the size of the
tested graph. If the graph fails the hyper-finiteness test,
then it is rejected. If it passes, then we know that by removing a
small proportion of the edges it is broken down to pieces of size $k$.
As these small pieces are still far from satisfying the property, we
can use a bounded number of random samples to actually find a subgraph
of the input that does not satisfy the property.

In Section \ref{s.proofs} we prove all of the stated result,
with the exception of Theorem~\ref{t.gap}, which is proved in
Section \ref{SecGap}. In particular, by combining a
theorem of Lipton and Tarjan with a result of Alon, Seymour and
Thomas \cite{AST}, regarding separators in minor-free graphs,
we get that $H$-minor free graphs are hyper-finite.
This facilitates a proof of Theorem~\ref{t1} from Theorem~\ref{t.m}.
In Section \ref{open} we discuss several open problems and conjectures for future
research.

\section{A Metric Criterion for Testability}\label{s.MC}

Let us slightly generalize the notion of property testing, and say
that two graph properties $A$ and $B$ are {\em distinguishable} if
there is a randomized algorithm that makes a bounded number of queries to
an input graph and accepts every graph in $A$ with probability at
least $2/3$ and rejects every graph in $B$ with probability at least
$2/3$. Then $\epsilon$-testing $A$ is equivalent to taking $B$
to be the set of graphs that are $\epsilon$-far from $A$.

Shortly, we will define a useful pseudo-metric\footnote{A
pseudometric on a set $X$ differs from a metric on $X$ in that a
pseudometric is allowed to be zero on pairs $(x,y)$ with $x \ne y$.
} on the set of graphs of bounded degree $d$. But first, we need to
set some terminology. A {\em rooted graph} is a pair $(G,v)$, where
$G$ is a graph and $v\in V(G)$. An {\em isomorphism between rooted
graphs} $(H,u)$ and $(G,v)$ is an isomorphism between the underlying
graphs $H$ and $G$ that maps $u$ to $v$. For a vertex $v\in V(G)$,
we denote by $B_G(v,r)$ the induced subgraph of $G$ whose vertices
consist of the vertices of $G$ at distance at most $r$ from $v$.
Consider a rooted graph $H$ and a finite graph $G$. Let $m^G_r(H)$
denote the number of vertices $v\in V(G)$ such that there is a
rooted-graph isomorphism from $\bigl(B_G(v,r),v\bigr)$ onto $H$. (Of
course, this is often zero.) Set $\mu^G_r(H):=m^G_r(H)/|V(G)|$, and
define a pseudometric $\dd_r$ by
$$
\dd_r(G,G'):=\sum_H \bigl| \mu^G_r(H)-\mu^{G'}_r(H) \bigr|\,,
$$
where the sum extends over all isomorphism types of rooted graphs.
Clearly, the number of terms that are nonzero is bounded by a
constant that depends only on $r$ and a bound on the degrees in $G$
and $G'$.
Observe that $\mu^G_r$ defines a probability measure on the set of rooted
graphs and that $\dd_r$ is monotone non-decreasing in $r$.
If $A$ and $B$ are graph families, we define
$\dd_r(A,B):=\inf\{\dd_r(G,G'):G\in A,G'\in B\}$.
The following proposition gives a metric condition for distinguishability.

\begin{proposition}\label{t.distinguish}
Let $A$ and $B$ be two graph properties having only graphs with
degrees at most $d$.
If there is some integer $R>0$  such that $\dd_R(A,B)>0$,
then $A$ and $B$ are distinguishable.
\end{proposition}

This gives a sufficient condition for $A$ to be testable:
if for every $\eps>0$ the set $B(\eps)$ of graphs that
are $\eps$-far from $A$ satisfies
$\sup_R \dd_R\bigl(A,B(\eps)\bigr)>0$, then $A$ is testable.

The converse to Proposition~\ref{t.distinguish} does
not hold. For example, if $A$ is the collection
of graphs with an even number of vertices and $B$ is the
collection with an odd number of vertices, then $A$ and $B$
are distinguishable in the current model in which the
number of vertices of the graph is given as the input.
However, it is not hard to come up with a natural model
for property testing in which the criterion given
by the proposition is necessary and sufficient.

The primary reason that makes hyper-finiteness so useful for property
testing is the following theorem.

\begin{theorem}\label{t.gap}
Fix $d,k\in\N_+$ and $\epsilon>0$. Let $A$ be the set of
$(\epsilon,k)$-hyper-finite graphs with vertex degrees
bounded by $d$, and let $B$ be the set of finite graphs with vertex
degrees bounded by $d$ that are not
$\bigl(4\,\epsilon\,\log(4\,d/\epsilon),k\bigr)$-hyper-finite.
Then there is some $R=R(d,k,\eps)\in\N_+$ such that $\dd_R(A,B)>0$.
\end{theorem}

Theorem \ref{t.gap} can actually be deduced from a recent result of
Schramm~\cite{Sc}, which deals with infinite unimodular graphs.
However, in the Section \ref{SecGap} we present a proof that is adapted to the finite
setting, gives quantitative bounds on $R$ and $\dd_R(A,B)$, and would
hopefully be more accessible.

\medskip

Let us now give an overview of the proof of Theorem \ref{t.gap}.
Assuming that $G$ is $(\epsilon,k)$-hyper-finite and that $G'$ is a
graph satisfying $\dd_R(G,G') \le \delta$ with some appropriate
$\delta>0$, we will show that $G'$ must be
$(4\,\epsilon\,\log(4\,d/\epsilon),k)$-hyper-finite. First, since
$G$ is $(\eps,k)$-hyper-finite, there is a subset $S$ of the edges
of $G$, of size at most $\epsilon |V(G)|$, that partitions $G$ into
connected components of size at most $k$. We then use the existence
of $S$ to construct a random $\tilde S\subset E(G)$, where each
connected component of $G\setminus\tilde S$ has at most $k$ vertices
and the expected size of $\tilde{S}$ is at most
$4\,\epsilon\,\log(3\,d/\epsilon)|V(G)|$.
 The most important feature of the way we will pick
$\tilde{S}$ is that it is local, in the sense that there is a finite bound
$R=R(d,\epsilon,k)$, such that the probability that
an edge $e$ is in $\tilde S$ only depends on the isomorphism
type of the pair $\bigl(B_G(e,R),e\bigr)$.
Now, if $G'$ is a graph satisfying
$\dd_R(G,G') \le \delta$, then locally it behaves almost exactly like
$G$ does. Therefore, choosing a set of edges $\tilde{S'}$ from $G'$
using the same process that was used for $G$, and removing it from
$G'$ should also partition $G'$ into connected components of size
at most $k$. Furthermore, the expected relative size of $\tilde{S'}$ in
$G'$ is close to that of $\tilde{S}$ in $G$, implying that $G'$ is
$(4\,\epsilon\,\log(4\,d/\epsilon),k)$-hyper-finite.

\medskip

We also have the following corollary which was mentioned in the
abstract.

\begin{corollary}\label{c.planardistance}
For every $\eps>0$ and $d\in\N_+$ there is an $R=R(\eps,d)\in\N_+$ such that
if $G$ and $G'$ are finite graphs with vertex degrees bounded by $d$,
$G$ is planar and $G'$ is $\eps$-far from being planar,
then $\dd_R(G,G')\ge 1/R$.
\end{corollary}

The proof of this corollary is an easy application of our arguments below,
and will be left to the reader.

\section{Most proofs}\label{s.proofs}

Here, we prove Proposition~\ref{t.distinguish} and Theorems~\ref{t.m} and~\ref{t1}.
The proof of Theorem~\ref{t.gap}, is postponed to Section \ref{SecGap}.

\proofof{Proposition \ref{t.distinguish}}
Suppose that for some integer $R>0$ and some positive $\delta>0$
we have $\dd_R(A,B)>\delta$.
 In that case, we can distinguish between $A$ and $B$
using the following algorithm. Let $\mathcal H$ denote the set of
(isomorphism types of) rooted graphs $(H,v)$ of radius $R$ around
$v$ and maximum degree at most $d$, and set $h:=|\mathcal H|$. Then,
given an input graph $G$ of size $n$, we estimate $\mu_R^G(H)$ for
every $H\in\mathcal H$, up to an additive error of $\delta/(2h)$,
with success probability $1-\frac{1}{4h}$. Here, we can apply an
additive Chernoff bound to deduce that to this end it is enough to
sample $O(\frac{h^2}{\delta^2}\log h)$ vertices, explore their
$R$-neighborhood, and compute the fraction of these vertices whose
neighborhood is isomorphic to $(H,v)$. Let $\hat{\mu}_R^G(H)$ be the
estimated values of $\mu_R^G(H)$, and recall that when computing
$\dd_R(G,G')$ we only need to consider rooted graphs $(H,v)$ of
radius $R$. The algorithm now checks\footnote{To do this, the
algorithm does not have to actually search all possible graphs. It
can just store a $\delta/4$-net of the set of all possible
$h$-tuples of values of $\mu^{G_A}(H)$, taken over all graphs $G_A
\in A$. This is a finite list.} if there exists a graph $G_A \in A$
for which $\sum_H \bigl| \hat{\mu}_R^G(H)-\mu^{G_A}_R(H) \bigr| \leq
\frac12\delta$. If that is the case, the algorithm declares that $G$
belongs to $A$, and otherwise it declares that $G$ belongs to $B$.
Note that since $\dd_R(A,B)>\delta$
and $\sum_H|\hat{\mu}_R^G(H)-\mu_R^G(H)| \leq \delta/2$
with high probability, the algorithm is unlikely to misclassify
graphs in $A$ or in $B$.
\QED

In order to deduce Theorem~\ref{t.m}, we will also need the following
lemma.

\begin{lemma}\label{l.mono}
Let $\mathcal P$ be a monotone hyper-finite graph property with degrees bounded
by $d$, let $\eps>0$, and let $B$ be the set of finite graphs with degrees
bounded by $d$ that are $\eps$-far from being in $\mathcal P$.  Then there
is a finite set of graphs $Z\subset\mathcal P$ and an $R\in\N_+$
such that
$\dd_R(B,\mathcal P\setminus Z)>0$.
\end{lemma}

\proof
Fix some $G'\in B$ and $G\in\mathcal P$.
Let $\eps_0\in(0,\eps/2)$ be sufficiently small so that
$4\,\epsilon_0\,\log(4\,d/\epsilon_0)<\eps/2$,
and let $k$ be such that each graph in $\mathcal P$ is
$(\eps_0,k)$-hyper-finite.
By Theorem~\ref{t.gap}, we have a lower
bound on $\dd_R(G,G')$ for some fixed $R=R(\mathcal P,\eps)$
if $G'$ is not $(\eps/2,k)$-hyper-finite.
So assume that $G'$ is $(\eps/2,k)$-hyper-finite.

Let $\mathcal K$ denote the set of all (isomorphism types of)
graphs $K$ such that each connected component of $K$
has at most $k$ vertices, and a disjoint union of finitely
many copies of $K$ is not
in $\mathcal P$.
Let $\mathcal K'$ denote the minimal graphs in $\mathcal K$;
that is, the set of graphs in $\mathcal K$ that do not
properly contain another graph of $\mathcal K$ as a subgraph.
Then the connected components of any $K\in\mathcal K'$ are
distinct, and hence $\mathcal K'$ is finite.

Assume, with no loss of generality, that the sizes of the
graphs in $\mathcal P$ are unbounded.
It then easily follows
that the graphs in $\mathcal K'$ do not have any connected
component with just one vertex.

Let $\mathcal H$ denote the (finite) collection of all rooted
connected graphs with
degrees bounded by $d$ such that each vertex is within distance at most
$k$ from the root vertex.
Given a graph $Y$, let $\mathcal H(Y)$ denote the rooted
graphs in $\mathcal H$ that contain $Y$ as a subgraph.
Suppose first that every $K\in \mathcal K'$ has a connected
component $Y$ such that
for every $H\in \mathcal H(Y)$ we have
\begin{equation}\label{e.smallfrequency}
\mu^{G'}_k(H)< \eps\,|\mathcal H|^{-1}\,d^{-k-1}/4\,.
\end{equation}
Then it follows that the same is true with $\mathcal K$ in place
of $\mathcal K'$.
Consequently, by removing all the edges of
$G'$ that are contained in rooted balls $H_v:=(B_{G'}(v,k),v)$ of $G'$ that
satisfy~\eqref{e.smallfrequency}
and furthermore removing a set of edges of size at most $\eps\,|V(G')|/2$
which cuts $G'$ to connected components of size at most $k$,
we arrive at a graph in $\mathcal P$.
By the choice of the right hand side of~\eqref{e.smallfrequency},
it follows that the number of edges removed from $G'$ is at most
$\eps\,|V(G')|/2$, which then contradicts the assumption that
$G'\in B$. Consequently, we find that there is some
$K\in\mathcal K'$ such that for every connected component
$Y$ of $K$ the inequality~\eqref{e.smallfrequency}
fails for some $H\in\mathcal H(Y)$.
For each such $Y$ let $H_Y$ be an $H\in\mathcal H(Y)$
for which~\eqref{e.smallfrequency} fails, and let
$\mathcal Q:=\{H_Y:Y\text{ is a connected component of }K\}$.

Suppose now that
$$
\dd_{k}(G,G')<\eps\,|\mathcal H|^{-1}\,d^{-k-1}/8\,.
$$
Then, in particular, we find that
$$
\mu^{G}_k(H)> \eps\,|\mathcal H|^{-1}\,d^{-k-1}/8
\qquad\qquad\forall
H\in\mathcal Q\,.
$$
Since $G\in\mathcal P$ and a disjoint union of finitely many copies
of $K$ is not in $\mathcal P$,
 this implies an
upper bound on $|V(G)|$, depending on $K$.
(In fact, depending on how many disjoint copies of $K$ are
sufficient to make a graph that is not in $\mathcal P$.)
However, we can maximize over the finitely many choices
of $K\in\mathcal K'$,
 to obtain an upper bound
that depends only on $\mathcal P$ and $\eps$.
Thus, except for a finite set of graphs $G\in\mathcal P$,
we have a lower bound on $\dd_k(G',G)$
or a lower bound on $\dd_R(G',G)$.
As $\dd_r$ is monotone non-decreasing in $r$, this proves the lemma.
\QED

\proofof{Theorem~\ref{t.m}}
Fix some $\eps>0$, and let $B$ be the set of graphs
that are $\eps$-far from $\mathcal P$.
Let $Z$ be the finite collection obtained from Lemma~\ref{l.mono}
with $B$ and $\mathcal P$ as above.
Then by Lemma~\ref{l.mono} and Proposition~\ref{t.distinguish},
$B$ and $\mathcal P\setminus Z$ are distinguishable.
Since $Z$ is finite,
the tester can easily test if the tested graph belongs to $Z$.
Thus, we get a tester for $\mathcal P$,
which first tests for membership in $Z$,
and then uses the tester which can tell apart $B$ from
$\mathcal P\setminus Z$. This completes the proof.
\QED

Recall the Lipton-Tarjan planar separator theorem~\cite{LT}, which
says that every planar graph $G$ has a set of vertices $V_0$ of size
$O(|V(G)|^{1/2})$, such that every connected component of
$G\setminus V_0$ has at most $2|V(G)|/3$ vertices. Lipton and Tarjan
used this theorem to show that the set of planar graphs with degree
bounded by $d$ is hyper-finite~\cite[Theorem 3]{LT1}. Alon, Seymour
and Thomas~\cite{AST} proved that the planar separator theorem as
stated above holds more generally in the class of $H$-minor free,
where $H$ can be any finite graph\footnote{The implicit constant in
the $O(\cdot)$ notation above is explicit in~\cite{LT}.
Also~\cite{AST} gives an estimate for this constant, which naturally
depends on $H$.}. The proof of~\cite[Theorem 3]{LT1} therefore gives
the following result.

\begin{proposition}\label{p.Hfree}
Fix $d\in\N_+$.
Let $H$ be a finite graph, and let $A$ be the set of graphs
that are $H$-minor free and have degrees bounded by $d$.
Then $A$ is hyper-finite.
\QED
\end{proposition}

\proofof{Theorem~\ref{t1}}
Since being $H$-minor free is a minor-closed property, it is enough to
prove the second claim of the theorem: that minor closed graph properties
are testable.
Let $\mathcal P$ be a minor closed graph property.
If $\mathcal P$ includes all graphs, then it is clearly testable.
Otherwise, suppose that $H\notin \mathcal P$.
Then $\mathcal P$ is $H$-minor free.
By Proposition~\ref{p.Hfree}, $\mathcal P$ is hyperfinite.
Since $\mathcal P$ is monotone,
the theorem now follows from  Theorem~\ref{t.m}.
\QED

\section{Proof of Theorem~\ref{t.gap}}\label{SecGap}

Let $G\in A$. Since $G$ is
$(\eps,k)$-hyper-finite, there is a set $S\subset E(G)$ such that
$|S|\le \eps\,|V(G)|$ and all the connected components of
$G\setminus S$ are of size at most $k$. Fix such a set $S$. With no
loss of generality, we assume that none of the edges of $S$ connects
two vertices from the same connected component of $G \setminus S$.
The proof strategy is to replace the set $S$ with a random set
$\tilde S$ that has no long-range dependencies and the probability
of an edge being in $\tilde S$ depends only on the local structure
of $G$ near the endpoints of the edge. It will then be easy to see
that a similar $\tilde S$ exists for $G_0$ if $\dd_R(G,G_0)$ is
small and $R$ is large.

\medskip

Let us first introduce some notation. If $K$ is a subgraph of a
graph $H$ and $K'$ is a subgraph of a graph $H'$, then we say that
the pair $(K,H)$ is isomorphic to the pair $(K',H')$ if there is an
isomorphism of $H$ onto $H'$ that takes $K$ onto $K'$. Isomorphisms
of triples $(H,K,J)$, $K,J\subset H$, are similarly defined. For a
graph $G$, we define $\mathcal K(G)$ as the set of vertex sets $K
\subset V(G)$ of size at most $k$ which span a connected graph in
$G$. Given the set $S$ that disconnects $G$ into connected
components of size at most $k$, we let $\mathcal K_S(G)$ consist of
those elements of $\mathcal K(G)$ that span a connected component in
$G\setminus S$. We stress that the elements in $\mathcal K(G)$ are only
required to span a connected subgraph in $G$, while the elements of
$\mathcal K_S(G)$ are required to span a {\em maximally} connected subgraph in
$G\setminus S$. In what follows it will be convenient to sometimes
identify an element $K \in K(G)$ with the (connected) graph spanned
by $K$.

Given a subgraph $H\subset G$ and $r\in\N$, let $N_r(H)$ denote the
subgraph of $G$ induced by the vertices at distance at most $r$ from
$H$. Let $N$ be some connected graph and let $K$ be a connected
subgraph of $N$ of size $|V(K)|\le k$. For $r\in\N_+$, let
$\Gamma_r(K,N)$ denote the set of $K'\in\mathcal K(G)$ such that
$\bigl(K',N_r(K')\bigr)$ is isomorphic to $(K,N)$. Set
$$
p_r(K,N):=
 \begin{cases}
 \frac{|\Gamma_r(K,N)\cap\mathcal K_S(G)|}{|\Gamma_r(K,N)|}
 &\text{if }\Gamma_r(K,N)\ne\emptyset\,, \\
 1  &\text{otherwise}\,.
 \end{cases}
$$
If $K\in\mathcal K(G)$, we abbreviate,
\begin{equation}
\label{e.abb}
p_r(K):= p_r\bigl(K,N_r(K)\bigr).
\end{equation}
For $v\in V(G)$, let $\mathcal K(v)$ denote the elements of
$\mathcal K(G)$ that contain $v$. Set
$$
q_r(v):=\sum_{K\in\mathcal K(v)} p_r(K)\,.
$$
Intuitively, $p_r(K)$ is an approximation
of the conditional probability that $K\in \mathcal K_S(G)$ given
the isomorphism type of the pair $\bigl(K,N_r(K)\bigr)$.
Since there is
exactly one $K\in\mathcal K(v)\cap\mathcal K_S(G)$,
we would expect that $q_r(v)$ is close to $1$ for most
vertices when $r$ is large. In the following,
we will need the fact that for some not too large
$R$, very few vertices have $q_R(v)$ small.
To this end, we have the following lemma.

\begin{lemma}\label{l.balance}
There exists some finite $R_1=R_1(d,\eps,k)$ and some
$R\in\N_+\cap[1,R_1]$ such that the number of vertices $v\in V(G)$
with $q_R(v)<1/2$ is at most $\eps\,|V(G)|/(2\,d)$.
\end{lemma}

We stress that $R$ itself may depend on $G$ and $S$,
but it is
bounded by $R_1$, which is not allowed to depend on $G$ and $S$. We
postpone the proof of the lemma and continue with the proof of
Theorem~\ref{t.gap}. Let $\mathcal K'$ denote a random subset of
$\mathcal K(G)$, where each $K\in\mathcal K(G)$ is in $\mathcal K'$
with probability $\min\bigl(2\,\log(2\,d/\eps)\,p_R(K),1\bigr)$,
independently, where $R$ is as provided by Lemma \ref{l.balance}.
Let
$$
S':=\bigcup_{K\in\mathcal K'}\partial K\,,
$$
where $\partial K$ denotes the set of edges of $G$ that are not in
$K$ but neighbor with some vertex in $K$. Let $W$ denote the set of
vertices that are contained in some $K\in\mathcal K'$, and $S''$
denote the set of all edges of $G$ having both vertices in
$V(G)\setminus W$. Finally, define $\tilde S:=S'\cup S''$. Clearly,
the connected components of $G\setminus\tilde S$ are of size at most
$k$. We now have to estimate $\E|\tilde S|$.

Let us start by estimating $\E|S''|$.
Consider an arbitrary vertex $v\in V(G)$ such that $q_R(v)\ge 1/2$.
Then
\begin{align*}
\Pb{v\notin W} &= \prod_{K\in \mathcal K(v)} \bigl(1-\Ps{K\in
\mathcal K'}\bigr)
\\&
\le \exp\Bigl(- \sum_{K\in\mathcal K(v)}
2\,\log(2\,d/\eps)\,p_R(K)\Bigr)
= \bigl(\eps/(2\,d)\bigr)^{2\,q_R(v)}.
\end{align*}
Therefore, if $q_R(v)\ge 1/2$, we have $\Pb{v\notin W}\le
\eps/(2\,d)$. On the other hand, the number of vertices satisfying
$q_R(v)<1/2$ is at most $\eps\,|V(G)|/(2\,d)$, by
Lemma~\ref{l.balance}. Thus, we have
$$
\E|V(G)\setminus W|\le
\frac{\eps|V(G)|}{2\,d}+\sum\Bigl\{\Ps{v\notin W}:v\in V(G),
q_R(v)\ge \frac 12\Bigr\} \le \frac{\eps|V(G)|}d\,.
$$
Since every edge in $S''$ is incident to a vertex in $V(G)\setminus
W$, we get $\E|S''|\le\eps\,|V(G)|$.

We now estimate $\E|S'|$.
For $K\in\mathcal K(G)$ set $L(K):= \Gamma_R\bigl(K,N_R(K)\bigr)$.
Observe that $K'\in L(K)$ if and only if $L(K')=L(K)$.
Let $K_1,K_2,\dots,K_m$ be a set of elements in $\mathcal K(G)$,
one from each equivalence class of the relation $L(K')=L(K)$.
Set $L_i:=L(K_i)$.
Note that $p_R(K)=p_R(K_i)$ for every $K\in L_i$,
and since $R>0$ also $|\p K|=|\p K_i|$ for every $K\in L_i$.
Therefore, the definition of $p_R$ gives
$$
\sum_{K\in L_i} p_R(K)\,|\partial K| = \sum_{K\in L_i}
\frac{|L_i\cap\mathcal K_S(G)|}{|L_i|}\,|\partial K_i| = \sum_{L_i\cap\mathcal K_S(G)} |\partial K|\,.
$$
By summing over $i$, we get,
$$
\sum_{K\in \mathcal K(G)} p_R(K)\,|\partial K|=\sum_i \sum_{K \in L_i}p_R(K)\,|\partial K|
=\sum_i \sum_{L_i \cap \mathcal K_S(G)} |\partial K|=\sum_{K \in \mathcal K_S(G)}|\partial K|=2|S|\,.
$$
As each element of $\mathcal K(G)$ is chosen with probability $\min\bigl(2\,\log(2\,d/\eps)\,p_R(K),1\bigr)$
to be in $\mathcal K'$,
we infer that
$$
\E|S'|\le 4\,\log(2\,d/\eps)\,|S| \le
4\,\eps\, \log(2\,d/\eps)\,|V(G)|
\,.
$$
Putting this together with our previous bound on $\E|S''|$, we obtain
\begin{equation}
\label{e.St}
\E|\tilde S| \le 4\,\eps\,\log(3\,d/\eps)\,|V(G)|\,.
\end{equation}

Now suppose that $G_0$ is any finite graph with degrees bounded by
$d$. We can define a random set of edges $S_0\subset E(G_0)$, as
follows. Let $\mathcal K'_0$ be a random set of elements of
$\mathcal K(G_0)$, where each $K\in \mathcal K(G_0)$ is placed in
$\mathcal K'_0$ with probability
$\min\bigl(2\,\log(2\,d/\eps)\,p_R(K,N_R(K)),1\bigr)$,
independently.
(Here, $N_R(K)$ refers to the neighborhood in $G_0$, of course.)
Define $S'_0:=\bigcup_{K\in \mathcal K'_0}\partial K$
and define $S''_0$ based on $G_0$ and $S'_0$ as $S''$ was defined
based on $G$ and $S'$.

In order to estimate $\E|\tilde S_0|$, we briefly consider the situation
in $G$ again.
Note that for every vertex $v \in V(G)$ the expected number of edges that
are incident with $v$ and belong to $\tilde S$ is completely determined
by its neighborhood of radius $r=R+k+1$, or, more precisely,
by the isomorphism type of the rooted graph $H(v):=\bigl(B(v,r),v\bigr)$.
Let $t_H$ be this expectation when $H$ is isomorphic to $H(v)$.
 Recall that
$\mu_r^G(H)\cdot |V(G)|$ is the number of vertices $v\in V(G)$
such that $H(v)$ is isomorphic to $H$.
 Then we can write
$$
\frac{2\E|\tilde S|}{|V(G)|}=\sum_{H}t_H \cdot \mu_r^G(H)\,.
$$
Similar considerations apply to $G_0$ and $\tilde S_0$, and hence
\begin{align*}
\frac{2\E|\tilde S_0|}{|V(G_0)|}=\sum_{H}t_H \cdot \mu_r^{G_0}(H) &= \sum_{H}t_H \cdot \mu_r^G(H) + \sum_{H}t_H \cdot (\mu_r^{G_0}(H)-\mu_r^G(H))\\
&= \frac{2\E|\tilde S|}{|V(G)|} + \sum_{H} t_H \cdot (\mu_r^{G_0}(H)-\mu_r^G(H))\\
& \leq  \frac{2\E|\tilde S|}{|V(G)|} + d \cdot \dd_{r}(G,G_0)\,,
\end{align*}
where we have used the fact that $t_H \leq d$ as $G_0$ is of bounded degree $d$.
Thus, $G_0\notin B$ if
\begin{equation}\label{DefDelta}
\dd_{r}(G,G_0)<(8/d)\,\eps\,\log(4/3)\;.
\end{equation}
This completes the proof. \QED

\proofof{Lemma \ref{l.balance}}
\def\hK{\mathcal Z}
\def\mG{\mathcal G}
Let $o$ be a uniformly random vertex chosen from $V(G)$.
Let $\mG$ be a set of finite rooted graphs, which has exactly
one representative for each isomorphism class of rooted graphs.
For $r\in\N$,
let $H_r$ denote the random element from $\mG$ that is
isomorphic to the rooted graph $\bigl(B(o,r),o\bigr)$.
Let $\hK$ denote the set of pairs $(H,K)$
such that $H\in\mG$ and $K$ is a connected subset of $V(H)$
of cardinality at most $k$ which contains the root of $H$.
It is important to note that unlike $\mG$, which does not
contain two isomorphic copies of the same graph, the set $\hK$ may contain two
isomorphic pairs $(H,K)$ and $(H,K')$.

For $Z=(H,K)\in\hK$, we
would like to define $\hat p_r(Z)$ as the conditional probability
given $H_r$ that $H=H_k$ and the isomorphic image of $K$ in
$B(o,k)$ is in $\mathcal K_S(G)$. But here is a slightly
annoying technical point. If $H_k$ has nontrivial
automorphisms, then there is no unique isomorphic image of $K$ in
$B(o,k)$. For this reason, given $o$ and $r$, we take a
random isomorphism $\gamma_r$, chosen uniformly among all
rooted-graph isomorphism from $H_r$ onto
$\bigl(B(o,r),o\bigr)$, with independent choices for each $r$. Let
$\beta_r$ denote the sequence of maps $(\gamma_{j+1}^{-1}\circ
\gamma_j:j=0,1,\dots,r-1)$.
Intuitively, $\beta_r$ tells us for each $r'<r$ how $H_{r'}$
sits inside $H_r$.

For $Z=(H,K)\in\hK$ let $A_Z$ denote the event that
$H_k=H$ and $\gamma_k(K)\in\mathcal K_S(G)$, and set
for every $Z\in\hK$
\begin{equation}
\label{e.pr}
\hat p_r(Z):= \Pb{A_Z\md H_r,\beta_r}.
\end{equation}
Clearly, exactly one of the events $A_Z$, $Z\in\hK$, holds.
Therefore,
\begin{equation}
\label{e.psum}
\sum_{Z\in\hK} \hat p_r(Z)=1
\end{equation}
holds for every $r\in\N$.
Since $(H_r,\beta_r)$ can be determined from $(H_{r+1},\beta_{r+1})$,
it is clear that
$$
\hat p_0(Z),~\hat p_{1}(Z),~\hat p_{2}(Z),\dots
$$
is a martingale for any fixed $Z\in\hK$.
(This is always the case when conditioning on finer and finer $\sigma$-fields.)
Since $\hat p_r(Z)\in[0,1]$ and $\Eb{\hat p_r(Z)}=\Pb{A_Z}=\hat p_0(Z)$,
we have $\Eb{\hat p_r(Z)^2}\le\hat p_0(Z)$.
Therefore, the variance of $\hat p_r(Z)$ is at most $\hat p_0(Z)$.

 Now recall that for a martingale $M_1,M_2,\dots$ we have the so called \lq\lq orthogonality for martingale differences\rq\rq, which says that
\begin{equation}
\label{e.md}
\sum_{j=1}^{n-1} \Eb{(M_{j+1}-M_{j})^2} = \Eb{M_n^2}-\Eb{M_1^2} \leq \Var(M_n).
\end{equation}
We apply this to the sequence $\hat p_{ki}(Z)$, $i=0,1,2,\dots$,
(which is also a martingale), sum over $Z\in\hK$
and use our above bound on the variance of $\hat p_r(Z)$, to obtain
for every $j\in\N_+$

\begin{equation}
\label{e.martdiff} \sum_{i=0}^{j-1} \sum_{Z\in \hK} \Eb{\bigl(\hat
p_{(i+1)k}(Z)-\hat p_{ik}(Z)\bigr)^2 } \le \sum_{Z\in \hK} \Var(\hat p_{jk}(Z)) \le\sum_{Z\in
\hK} \hat p_0(Z)=1\,.
\end{equation}

It is now time to make the connection between $\hat p_r$
and $p_r$.
Fix some $i\in\N_+$ and $Z=(H,K)\in \hK$.
Suppose that $H_k=H$. Then
$$
\gamma_{(i+1)k}(H_{(i+1)k})=B\bigl(o,(i+1)k\bigr)\supset N_{ik}\bigl(\gamma_k(K)\bigr)
\supset
\gamma_{ik}(H_{ik})=B(o,ik)\,,
$$
where, as before $N_r(U)$ denotes the $r$-neighborhood of $U$.
Let $Y$ denote the isomorphism type of the triple
$\bigl(N_{ki}\bigl(\gamma_k(K)\bigr),\gamma_k(K),o\bigr)$. Then
clearly $Y$ can be determined from $Z$, $H_{(i+1)k}$ and
$\beta_{(i+1)k}$. Moreover, $Y$ determines the isomorphism type of
$\bigl(H_{ik},\gamma_{ik}^{-1}\circ\gamma_k(K),o_{ik}\bigr)$ and
therefore determines $\hat p_{ik}(Z)$. Consequently, the three term
sequence
\begin{equation}
\label{e.m2}
\hat p_{ik}(Z), ~1_{\{H_k=H\}}\,\Pb{A_Z \md Y}, ~\hat p_{(i+1)k}(Z)
\end{equation}
is a martingale.
We abbreviate  $\bar p_{ik}(Z):=1_{\{H_k=H\}}\,\Pb{A_Z\md Y}$.
It is easy to see that
$$
\bar p_{ik}(Z)=\begin{cases} p_{ik}\bigl(\gamma_k(K)\bigr) & H=H_k\,,\\
 0 & H\ne H_k
\,.
\end{cases}
$$
Therefore,
\begin{equation}
\label{e.q}
q_{ik}(o)=\sum_{Z\in\hK} \bar p_{ik}(Z)\,.
\end{equation}
Since the three random variables $\hat p_{ik}(Z), ~\bar p_{ik}(Z), ~\hat p_{(i+1)k}(Z)$ form a martingale,
we get from~\eqref{e.md} that
\begin{equation}
\label{e.near}
\begin{aligned}
\EB{\bigl(\hat p_{(i+1)k}(Z)-\hat p_{ik}(Z)\bigr)^2}
&
 = \EB{\bigl(\hat p_{(i+1)k}(Z)-
\bar p_{ik}(Z)\bigr)^2}+
\EB{\bigl(\bar p_{ik}(Z)-
\hat p_{ik}(Z)\bigr)^2}
\\ &
\ge
\EB{\bigl(\bar p_{ik}(Z)-\hat p_{ik}(Z)\bigr)^2}   .
\end{aligned}
\end{equation}
Recall the definition of $\mathcal K(v)$ from just below~\eqref{e.abb}.
There is clearly a finite upper bound\footnote{One way of seeing this
is to observe that every tree of size $k$ can be realized as a closed walk of length $2k$; just
double each edge and take an Euler tour.
Therefore, the number of sub-trees of size at most $k$ that contain the root
of a graph, whose maximum degree is $d$, is at most $d^{2k}$.} $t=t(k,d) \leq d^{2k}$ on $|\mathcal K(v)|$.
Observe that
$$
\bigl|\{(H,K)\in\hK: H_k=H\}\bigr|=\bigl|\mathcal K(o)\bigr|\le t\,.
$$
Thus, there are at most $t$ different $Z\in\hK$ for which
$\bar p_{ik}(Z)-\hat p_{ik}(Z)\ne 0$.
Cauchy-Schwarz therefore gives
\begin{equation}
t\sum_{Z\in\hK}\Bigl(\bar p_{ik}(Z)-\hat p_{ik}(Z)\Bigr)^2
\ge
\Bigl(\sum_{Z\in\hK}\bar p_{ik}(Z)-\sum_{Z\in\hK}\hat p_{ik}(Z)\Bigr)^2.
\end{equation}
Taking expectation on both sides and using linearity of expectation we get
\begin{equation}
\label{e.CS}
t\sum_{Z\in\hK}\EB{\Bigl(\bar p_{ik}(Z)-\hat p_{ik}(Z)\Bigr)^2}
\ge
\EB{\Bigl(\sum_{Z\in\hK}\bar p_{ik}(Z)-\sum_{Z\in\hK}\hat p_{ik}(Z)\Bigr)^2}.
\end{equation}
We now sum~\eqref{e.near} over all $Z\in\hK$ and apply~\eqref{e.CS} and
then
\eqref{e.psum} and~\eqref{e.q}, to get
$$
t\,\sum_{Z\in\hK}
\EB{\bigl(\hat p_{(i+1)k}(Z)-\hat p_{ik}(Z)\bigr)^2}
\ge
\EB{
\Bigl(\sum_{Z\in\hK}\bar p_{ik}(Z)-\sum_{Z\in\hK}\hat p_{ik}(Z)\Bigr)^2}
=\Eb{\bigl(q_{ik}(o)-1\bigr)^2}.
$$
Now sum over $i$ and apply~\eqref{e.martdiff},  to obtain
$$
t \ge
\sum_{i=1}^{j-1}
\EB{\bigl(q_{ik}(o)-1\bigr)^2}
=|V(G)|^{-1}\sum_{v\in V(G)}\sum_{i=1}^{j-1}
\bigl(q_{ik}(v)-1\bigr)^2\,.
$$
To complete the proof we can take
\begin{equation}\label{DefR}
j=k+8\,d\,t\,k/\eps \leq 10\,k\,d^{2k+1}/\epsilon \;,
\end{equation}
and the same choice works for $R_1$ in the statement of the lemma.
\QED

\section{Concluding Remarks and Open Problems}\label{open}

\begin{enumerate}

\item Let us give an explicit upper bound for the
dependence on $\epsilon$ of the number of queries required
to test the property of being $H$-minor free for any fixed
connected graph $H$. We assume here that the maximum degree $d$ is a
constant and $H$ is fixed. Given a graph $G$ and $\epsilon > 0$, let us choose
$\epsilon_0$ and $k$ as in the proof of Lemma \ref{l.mono}. It is
clear that $1/\epsilon_0=\mathrm{poly}(1/\epsilon)$ and the result of
\cite{AST} together with \cite[Theorem 3]{LT1} give that
$k=\mathrm{poly}(1/\epsilon)$. By Theorem~\ref{t.gap} and
its proof we know that
if $G$ is not $(\epsilon/2,k)$-hyper-finite then for every $G'$ that
is $H$-minor-free we have $\dd_R(G,G')=\Omega(\epsilon)$ for some
fixed $R=2^{O(k)}/\epsilon=2^{\mathrm{poly}(1/\epsilon)}$. Here we
are using the bounds that appear in (\ref{DefDelta}) and
(\ref{DefR}). The algorithm first tries to distinguish between the
case that the input is $(\epsilon_0,k)$-hyper-finite from the case
that it is not $(\epsilon/2,k)$-hyper-finite. By
the proof of Proposition
\ref{t.distinguish} this can be done with
$\mathrm{poly}(h/\epsilon)$ queries, where $h$ is the number of
graphs of bounded degree $d$ and radius $R$. Clearly
$h=2^{2^{O(R)}}=2^{2^{2^{\mathrm{poly}(1/\epsilon)}}}$. If the
algorithm finds that $G$ is not $(\epsilon_0,k)$-hyper-finite it
rejects $G$. Otherwise, the algorithm samples $\mathrm{poly}(1/\epsilon)$
vertices, and for each vertex it explores its neighborhood of radius
$k$. This requires $2^{O(k)}=2^{\mathrm{poly}(1/\epsilon)}$ queries.
If any of the neighborhoods explored is not $H$-minor free the
algorithm rejects, otherwise it accepts. The algorithm clearly
accepts $G$ with high probability if it is $H$-minor free, as such a
$G$ is $(\epsilon_0,k)$-hyper-finite. If $G$ is $\epsilon$-far from
being $H$-minor free and not $(\epsilon/2,k)$-hyper-finite then it
is rejected with high probability in the first step. If not then it
is easy to see that $\Omega(\epsilon n)$ of its vertices belong to
connected subgraphs of $G$ of radius at most $k$ that are not
$H$-minor free. Hence, $G$ is rejected in the second step.

If $H$ is not connected, then after testing hyperfiniteness,
we can test for being $H'$-minor free for each connected
component $H'$ of $H$. It is then easy to determine
if $G$ is $H$-minor free or not. Similar reasoning applies
in the case of a property determined by a finite list of
forbidden minors. By the Graph-Minor Theorem~\cite{RobSey2},
this includes all minor closed properties. Thus,
every minor closed property can be tested using
 $2^{2^{2^{\mathrm{poly}(1/\epsilon)}}}$ queries.

\item
Another model of
property testing is when the number of edges is arbitrary, and the
error is relative to the number of edges \cite{PR}; that is, a graph
with $m$ edges is $\epsilon$-far from ${\cal P}$ if we have to
modify $\epsilon m$ of its edges to get a graph satisfying ${\cal
P}$. This raises the following:

\begin{problem}
What is the query complexity of testing minor-closed properties in
graphs with bounded {\em average} degree?
\end{problem}

It is clear that $o(n^{1/2})$ queries are not enough,
for they will not suffice to detect a clique on $|V(G)|^{1/2}$
of the vertices of the graph.

\item As we have argued in Subsection \ref{compare}, minor closed properties
are trivially testable in the dense graph model using
$O(1/\epsilon)$ queries even with one-sided error. The intuition is
that if a hereditary graph property is testable, then the reason
must be that one can find a ``proof'' that a graph does not satisfy
the property, in the form of a subgraph that does not satisfy it
(which implies that the graph does not satisfy it). This intuition
turns out to be correct in dense graphs \cite{ASHer}. However, it is
easy to see that in the bounded degree model, we cannot expect to
have such proofs. For example, if $G$ is a bounded degree expander
of girth $\Omega(\log n)$, then for any fixed $H$ with at least one
cycle, we have that $G$ is $\Omega(1)$-far from being $H$-minor
free\footnote{This follows, for example, from the separator theorem
of \cite{AST} for $H$-minor free graphs, mentioned earlier.} when it
is sufficiently large, but on the other hand, every subgraph of $G$
of size $o(\log n)$ is a tree, and in particular $H$-minor free.
Thus if $H$ is not a tree, then one cannot test $H$-minor freeness
with $o(\log n)$ queries and one-sided error. In fact, a much
stronger $\Omega(\sqrt{n})$ lower bound can be deduced by adapting
an argument from \cite{GR}. We raise the following conjecture, stating
that the $\Omega(\sqrt{n})$ lower bound is tight.

\begin{conjecture}\label{Conj1}
For every $H$, being $H$-minor free can be tested
in the bounded degree setting with one-sided
error and query complexity $\tilde{O}(\sqrt{n})$.
\end{conjecture}

If the conjecture is true, then the Graph Minor Theorem~\cite{RobSey2}
implies that the same is true for any minor-closed graph property.

\item As we have briefly mentioned earlier, it does not seem like
our approach can lead to testing algorithms for the property of
being $H$-minor free whose query complexity is polynomial in
$1/\epsilon$. It seems interesting to further investigate the
possibility of coming up with such a tester.

\item Our main result here is that we can distinguish in constant
time between graphs satisfying a minor closed property, from those
that are far from satisfying it. Can we also estimate in constant
time the fraction of edges that need to be removed in order to make
the graph satisfy the property?

\item It would be interesting to obtain an explicit description of
the difference between the frequencies  of local neighborhoods
that one sees in $(\eps,k)$-hyper-finite graphs versus
graphs that are not $(\eps',k)$-hyper-finite, where $\eps\ll\eps'$.
\end{enumerate}

\paragraph{Acknowledgements:} We would like to thank David Wilson and
G\'abor Elek for helpful discussions.


\begin{thebibliography}{99}


\bibitem{AFNS}
N. Alon, E. Fischer, I. Newman and A. Shapira, A combinatorial
characterization of the testable graph properties: it's all about
regularity, Proc. of STOC'06, 251--260. Also, SIAM J. on Computing,
to appear.

\bibitem{AST}
N. Alon, P. D. Seymour and R. Thomas, A separator theorem for
non-planar graphs, J. Amer. Math. Soc. 3 (1990), 801--808.

\bibitem{ASHer}
N. Alon and A. Shapira, A characterization of the (natural) graph
properties testable with one-sided error, Proc. of FOCS 2005,
429--438. Also, SIAM J. on Computing, to appear.

\bibitem{ASSurvey}
N. Alon\ and\ A. Shapira, Homomorphisms in graph property testing,
in {\it Topics in discrete mathematics}, 281--313, Springer, Berlin,
2006.

\bibitem{ASSep}
N. Alon and A. Shapira, A separation theorem in property testing,
Combinatorica, to appear.

\bibitem{BLR}
M. Blum, M. Luby and R. Rubinfeld, Self-testing/correcting with
applications to numerical problems, JCSS 47 (1993), 549--595.

\bibitem{Bod}
H. L. Bodlaender, Discovering treewidth, Lecture Notes in Computer
Science, 3381 (2005), 1--16.

\bibitem{BOT}
A.~Bogdanov, K.~Obata, and L.~Trevisan, A lower bound for testing
3-colorability in bounded-degree graphs, FOCS 2002, 93--102.

\bibitem{BCLSSV06}
C.~Borgs, J.~Chayes, L.~Lov\'asz, V.~T.~Sos, B.~Szegedy, and
K.~Vesztergombi, Graph limits and parameter testing,
STOC 2006, 261--270.

\bibitem{CSS}
A. Czumaj, A. Shapira, and C. Sohler, Testing hereditary properties
of non-expanding bounded-degree graphs, submitted (full version of
\cite{CS1}).

\bibitem{CzSo}
A. Czumaj and C. Sohler, Testing expansion in bounded-degree graphs,
Proc. of FOCS 2007, 570-578.

\bibitem{CzSo2}
A. Czumaj and C. Sohler, Sublinear-time algorithms, Bulletin of the
EATCS, 89 (2006) 23--47.

\bibitem{CS1}
A. Czumaj and C. Sohler, On testable properties in bounded degree
graphs, Proc. of SODA 2007, 494--501.

\bibitem{DHK}
E. Demaine, M. Hajiaghayi and K. Kawarabayashi, Algorithmic graph
minor theory: decomposition, approximation, and coloring, Proc. of
FOCS 2005, 637-646.

\bibitem{Dies}
R. Diestel, {\bf Graph Theory} (Third Edition), Springer, Heidenberg, 2005.

\bibitem{DST}
C. Domingo and J. Shawe-Taylor, The Complexity of Learning
Minor-Closed Graph Classes, Proc. of Algorithmic Learning Theory
1995, 249--260.

\bibitem{DF}
R. Downey and M. Fellows, {\bf Parameterized Complexity}, Springer,
1999.

\bibitem{Elek}
G. Elek, The combinatorial cost, to appear in L'Enseignement Math\'ematique

\bibitem{Elek2}
G. Elek, A Regularity lemma for bounded degree graphs and its
applications: parameter testing and infinite volume limits,
arXiv:07112800, 2007.

\bibitem{F}
E. Fischer, The art of uninformed decisions: A primer to property
testing, The Computational Complexity Column of The Bulletin of
the European Association for Theoretical Computer Science 75
(2001), 97--126.

\bibitem{GGR}
O.~Goldreich, S.~Goldwasser and D.~Ron, Property testing and its
connection to learning and approximation, J. of the ACM 45(4):
653--750 (1998).

\bibitem{GR}
O.~Goldreich and D.~Ron, Property Testing in Bounded-Degree Graphs,
Algorithmica 32 (2002), 302--343.

\bibitem{GR99}
O.~Goldreich and D.~Ron, A sublinear bipartiteness tester for
bounded degree graphs, Combinatorica, 19 (1999), 335--373.


\bibitem{GT}
O. Goldreich and L. Trevisan, Three theorems regarding testing graph
properties, Random Structures and Algorithms, 23(1):23--57, 2003.


\bibitem{HT}
J. E. Hopcroft and R. E. Tarjan, Efficient planarity testing, J. of
the ACM, 21 (1974), 549--568.

\bibitem{K1}
A. Kostochka, The minimum Hadwiger number for graphs with a given
mean degree of vertices, Metody Diskret. Analliz. 38 (1982), 37--58
[in Russion].

\bibitem{K2}
A. Kostochka, A lower bound for the Hadwiger number of graphs by
their average degree, Combinatorica, 4 (1984), 307--316.

\bibitem{Kur}
K. Kuratowski, Sur le probl\`{e}me des courbes gauches en topologie,
Fund. Math. 15: 271--283.


\bibitem{LT}
R. J. Lipton and R. E. Tarjan, A separator theorem for planar graphs,
SIAM J. Appl. Math., 36 (1979), 177--189.

\bibitem{LT1}
R. J. Lipton and R. E. Tarjan, Applications of a planar separator
theorem, SIAM J. on Computing 9 (1980), 615--627.


\bibitem{L}
L. Lov\'asz, Graph minor theory, Bull. Amer. Math. Soc. 43 (2006),
no. 1, 75--86.

\bibitem{PR}
M. Parnas and D. Ron, Testing the diameter of graphs, Random
structures and algorithms, 20 (2002), 165--183.

\bibitem{RobSey1}
N. Robertson and P.D.~Seymour: Graph minors XIII. The disjoint paths problem, J. Combin.
Theory Ser. B 63 (1995), 65--110.

\bibitem{RobSey2}
N. Robertson and P. D. Seymour, Graph minors. XX. Wagner's
Conjecture, J. Combin. Theory Ser. B 92 (2004), 325--357.

\bibitem{Rub}
R. Rubinfeld, Sublinear time algorithms, in {\it International Congress of Mathematicians. Vol. III},
1095--1110, Eur. Math. Soc., Z\"urich, 2006.

\bibitem{Ron} D. Ron,
Property testing, in: P. M. Pardalos, S. Rajasekaran, J. Reif and
J. D. P. Rolim, editors, {\em Handbook of Randomized Computing},
Vol. II, Kluwer Academic Publishers, 2001, 597--649.

\bibitem{RS} R.~Rubinfeld and M.~Sudan,
Robust characterization of polynomials with applications to program
testing, SIAM J. on Computing, 25 (1996), 252--271.

\bibitem{Sc}
O.~Schramm, Graph sequences with hyperfinite limits are hyperfinite, arXiv:0711.3808, 2007.

\bibitem{Sz} E.~Szemer\'edi,
Regular partitions of graphs, In: {\em Proc.\ Colloque Inter.\
CNRS} (J.~C.~Bermond, J.~C.~Fournier, M.~Las~Vergnas and
D.~Sotteau, eds.), 1978, 399--401.

\bibitem{T1}
A. Thomason, An extremal function for contractions of graphs, Math.
Proc. Cambridge Philos. Soc. 95 (1984), 261--265.

\bibitem{T2}
A. Thomason, The extremal function for complete minors, J. of
Combinatorial Theory Series B, 81 (2001), 318--338.

\bibitem{W}
K. Wagner, \"{U}ber eine eigenschaft der ebenen komplexe, Math. Ann. 114 (1937), 570--590.

\end{thebibliography}
\end{document}